\documentclass[12pt,twoside]{article}
\long\def\symbolfootnote[#1]#2{\begingroup%
\def\thefootnote{\fnsymbol{footnote}}\footnote[#1]{#2}\endgroup}

\usepackage{amsfonts, amssymb, amsmath, amsthm}
\DeclareMathOperator*{\esssup}{ess\,sup}
\usepackage[colorlinks=true, pdfstartview=FitV, linkcolor=blue,
            citecolor=blue, urlcolor=blue]{hyperref}
\usepackage[usenames]{color}
\definecolor{Red}{rgb}{0.7,0,0.1}
\definecolor{Green}{rgb}{0,0.7,0}
\usepackage{accents}
\usepackage{comment}

\newcommand{\dom}{\mathcal M}

\title
{Local Existence of Strong Solutions to the $3D$  Zakharov-Kuznestov Equation in a Bounded Domain}
\author{Chuntian Wang 
}

\date{}

\numberwithin{equation}{section}

\newtheorem{thm}{Theorem}[section]
\newtheorem{lem}{Lemma}[section]
\newtheorem{prop}{Proposition}[section]

\newtheorem{rem}{Remark}[section]

\begin{document}


\maketitle

\vskip-4mm

\centerline{\footnotesize{\it  Department of Mathematics and The Institute for Scientific Computing and Applied Mathematics    }}
\vskip-1mm
\centerline{\footnotesize{\it Indiana University, Bloomington, IN 47405}}

\vskip 0mm
\centerline{\footnotesize{\it \,\, email: \url{wang211@umail.iu.edu}}}

%
%
%
%
%

\begin{center}
\large
\date{\today}
\end{center}

\vskip4mm

\tableofcontents

\newpage

\begin{abstract}
 We consider here the local existence of strong solutions for the Zakharov-Kuznestov (ZK) equation posed in a limited domain $\dom=(0,1)_{x}\times(-\pi /2, \pi /2)^d,$ $ d=1,2$. We  prove that in space dimensions $2$ and $3$, there  exists a strong solution on a short time interval, whose length only depends on the given data.
We use the  parabolic regularization of the ZK equation as in \cite{SautTemamChuntian} to  derive the global and local bounds independent of $\epsilon$ for various norms of the solution. In particular, we derive the local bound of the nonlinear term  by a singular perturbation argument. Then we can  pass to the limit and hence deduce  the local existence of strong solutions.
\end{abstract}
{\noindent \small
 {\it \bf Keywords: Zakharov-Kuznetsov equation, Korteweg-de Vries equation, Local existence of  strong solutions in $3D$\\}


\section{Introduction}
\label{sec:introduction}
The Zakharov-Kuznetsov equation
\begin{equation}\label{eq0r}
\dfrac{\partial u}{\partial t} + \Delta \dfrac{\partial  u}{ \partial x  } +  c \dfrac{  \partial u  }{ \partial x  }+ u \dfrac{  \partial u  }{ \partial x  }   =f,
\end{equation}
where $u=u(x,  x^\perp, t)$,  $x^\perp =y$ or $x^\perp=(y,z)$, describes the propagation of nonlinear ionic-sonic waves in a plasma submitted to a magnetic field directed along the $x$-axis.  Here $c>0$ is the sound velocity.  It has been derived formally in a long wave, weakly nonlinear regime from  the Euler-Poisson system in \cite{ZK} and
 \cite {LSp}. A rigorous derivation is provided in \cite{LLS}.  For more general physical references, see   \cite{bonaparitchard1} and \cite{bonapritchatd2}.
When $u$ depends only on $x$ and $t$, (\ref{eq0r}) reduces to the classical Korteweg-de Vries (KdV) equation.



%

Concerning the initial and  boundary value problems of  the Korteweg-de Vries equation posed on a bounded interval $(0,L)$, we refer the interested readers to e.g. \cite{bonasunzhang}, \cite{colinghida} and \cite{colingisclon}.

The initial and boundary value problem associated  with (\ref{eq0r}) has been studied
 in the half space $\{ (x, y):\, x>0 \}$ (\cite{Fam2}), on a strip like $\{(x, y):\, x \in \mathbb R, \, 0<y <L \}$ (\cite{BaykovaFaminskii}) or $\{(x, x^\perp):\, 0<x<1, \,x^\perp \in \mathbb R^d, \, d=1,2 \}$ (\cite{Fam3} and \cite{SautTemam}),
  and in a rectangle $ \{ (x, x^\perp): \, 0<x<1, \, x^\perp \in (-\pi/2, \pi/2)^d, \, d = 1, 2 \}$  (\cite{SautTemamChuntian}). Specifically in \cite{SautTemamChuntian},  the authors  have established, for  arbitrary large initial data, the existence of global weak solutions in space dimensions  $2$ and $3$ ($d=1$ and $ 2$ respectively) and a result of uniqueness of such solutions in the two-dimensional case.

As for the existence of  strong solutions,   the global existence  in space dimension $2$ has been proven in a half strip $\{ (x, y):\, x>0,\, y \in (0, L) \}$ in \cite{LT}.  The existence and  exponential decay of regular solutions to the linearized ZK equation in a rectangle $ \{ (x, y): \, x \in (0, L), \, y \in (0, B) \}$ has been studied in \cite{DoroninLarkin}.

However, to the best of our knowledge, there has been no result so far     for the local existence of strong solutions in $3D$ in a limited domain. In the present article  we prove  the short time existence of strong solutions in a $3D$   rectangular domain as in \cite{SautTemamChuntian} .

%
%

The  article is organized as follows.

Firstly we introduce the basic settings of the equation and the related functional spaces in Section \ref{sec0}.

 Secondly we introduce the  parabolic regularization  as in \cite {SautTemam} and \cite{SautTemamChuntian}. We derive the global bounds on   $u^\epsilon$ independent of $\epsilon$ in Section \ref{secglobal}. Then we derive the local bounds on $[0, T_{\ast})$, where  $T_{\ast} >0$ depends only on the data (Sections \ref{seclocal} and  \ref{sing}). In particular, we use a singular perturbation argument in the $x$ direction to deduce the local bounds on the nonlinear term in Section \ref{sing}.

 Finally we pass to the limit on the regularized equation
 and obtain the local existence of a strong solution (Section \ref{secexistence}).


%
%


In the Appendix, we recall from \cite{SautTemam} and \cite{SautTemamChuntian} a trace result concerning a singular perturbation problem in the $x$ direction.

%
%
%
%
%
%
%


\section{The ZK equation in a rectangle}
\label{sec0}
We aim to  study the ZK equation:
\begin{equation}\label{eq0}
\dfrac{\partial u}{\partial t} + \Delta \dfrac{\partial  u}{ \partial x  } +  c \dfrac{  \partial u  }{ \partial x  }+ u \dfrac{  \partial u  }{ \partial x  }   =f,
\end{equation}
in a rectangle or  parallelepiped domain  in $\mathbb R^n$  with $n= 2$ or $ 3$, denoted as  $\dom=(0,1)_{x}\times(-\pi /2, \pi /2)^d,$ with $ d=1 $ or $2$, $\Delta u=u_{xx}+\Delta^\perp u$, $\Delta^\perp u=u_{yy}$ or $u_{yy}+u_{zz}$.
In the sequel we will use the notations $I_{x}=(0,1)_{x},$ $ I_{y}=(- \pi /2,  \pi /2)_{y}$,  $I_{z}=(- \pi /2, \pi /2)_{z}$, and $I_{x^\perp} =I_y$ or $ I_{y} \times I_{z}$. We assume the boundary conditions on $x=0, 1$
 to be
\begin{equation}\label{eq3}
 u(0, \, x^\perp,\, t)=u(1, \, x^\perp,\, t)=u_x( 1, \, x^\perp,\, t )=0,
\end{equation}
and the initial condition reads:
\begin{equation}\label{eq33}
 u(x, \, x^\perp,\, 0)=u_0(x, \, x^\perp).
\end{equation}
We also need  suitable boundary conditions in the $y$ and $z$ directions.  As in \cite{SautTemamChuntian},  we will choose  either the Dirichlet boundary conditions
\begin{equation}\label{eq107}
  u=0 \mbox{ at } y=\pm\frac{\pi}{2}\,\,\,(z=\pm\frac{\pi}{2} ),
\end{equation}
or the periodic boundary conditions
\begin{equation}\label{eq1-26p}
 u\big|^ {y= \frac {\pi}{2}}_{y=-\frac {\pi}{2}}=u_y\big|^ {y= \frac {\pi}{2}}_{y=-\frac {\pi}{2}}=0\,\,\,(u\big|^{z= \frac {\pi}{2}} _{z=-\frac {\pi}{2}}=u_z\big|^{z= \frac {\pi}{2}} _{z=-\frac {\pi}{2}}=0).
\end{equation}
 We will study the initial and boundary value problem (\ref{eq0})-(\ref{eq33}) supplemented with the boundary condition (\ref{eq107}), that is, in the  Dirichlet case,  and we will make some remarks on the extension to the periodic boundary condition case.

We will denote by $|\cdot|$ and $(\cdot, \cdot) $ the norm and the inner product of $L^2 (\dom)$.

We use the following functional space in the sequel:
\begin{align*}
A u: = \Delta u_ x + cu_x,\quad\quad \,\,\,\, &\forall\,\, u\in D(A),
\end{align*}
where  $D(A)=\left \{u\in L^2 (\dom): \,\,Au\in L^2 (\dom), \,\, u=0 \,\,\mbox{on}\,\, \partial \dom, \,\, u_x( 1, \, x^\perp,\, t )=0  \right\}$. Note that the trace theorem proven in \cite{SautTemamChuntian} shows that if $u\in L^2 (\dom)$ and $Au \in L^2 (\dom)$ then the traces of $u$ on $\partial \dom$, and of  $u_x$ at $x=1$ make sense.

We also consider the space
 \begin{equation}\label{xi1}
\Xi = \left \{ u \in H^2(\dom) \cap H_0^1(\dom), \,\,u_x  \big|_{x=1}=0   \right  \},
\end{equation}
and endow this space  with the scalar  product and  norm  $[ \cdot, \cdot ]_2$ and $ [ \cdot ]_2$,
\begin{align*}
[ u  , v ]_2 =  (u_{xx}, v_{xx}) + (u_{yy}, v_{yy}) + (u_{zz}, v_{zz}),
\end{align*}
\begin{equation}\label{[]}
[ u ]_2^2  =    \left| u_{xx} \right|^2 + \left| u_{yy} \right|^2 +\left| u_{zz} \right|^2  ,
\end{equation}
which make it a Hilbert space, thanks to the Dirichlet boundary condition and the     Poincar\'e inequality and elliptic boundary regularity (see \cite{SautTemam} and \cite{SautTemamChuntian}). Furthermore by a result proven in \cite{SautTemamChuntian}, we know that $D(A) \subset H^2 (\dom) \cap H_0^1(\dom)$, and hence
\begin{equation}\label{da2}
D(A) \subset \Xi.
\end{equation}

\section{Parabolic regularization}
\label{secparabolic}
For the sake of simplicity we will only treat the more complicated case when $d=2$; the case when $d=1$ is easier.
To begin with, we recall the  parabolic regularization introduced in \cite{SautTemam} and  \cite{SautTemamChuntian},
that is, for $\epsilon >0$ ``small'', we consider the parabolic equation,
\begin{eqnarray}\label{10-11}\displaystyle
\begin{cases}
&\,\,\, \dfrac{\partial u^\epsilon}{\partial t}+ \Delta \dfrac{\partial u^\epsilon }{\partial x}+ c \dfrac{\partial u^\epsilon }{\partial x} +     u^\epsilon \dfrac{\partial u^\epsilon }{\partial x}  + \epsilon \, L u^\epsilon      =f,\\
&\,\,\,u^\epsilon (0)=u_0,
\end{cases}
\end{eqnarray}
where
  \begin{equation*}
L u^\epsilon : =  \dfrac{\partial^4 u^\epsilon}{\partial x^4}   +
 \dfrac{\partial^4 u^\epsilon}{\partial y^4}   + \dfrac{\partial^4 u^\epsilon}{\partial z^4} ,
\end{equation*}
  supplemented with the boundary conditions (\ref{eq3}), (\ref{eq107}) and the additional boundary conditions
\begin{equation}\label{10-4}
u^\epsilon _{yy}\big|_{ {y=\pm \frac {\pi}{2}}  }=u^\epsilon _{zz}\big|_{ {z=\pm \frac {\pi}{2}}  }=0,
\end{equation}
\begin{equation}\label{10-2}
u^\epsilon _{xx}\big|_{x=0}  = 0.
\end{equation}
Note that since $u^\epsilon_{yy}= u^\epsilon_{zz}=0$ at $x=0$,  (\ref{10-2}) is equivalent to
\begin{equation*}
\Delta u^\epsilon \big|_{x=0}  = 0.
\end{equation*}

It is a classical result (see e.g. \cite{JLLions},     \cite{MR0241822} or also  \cite{SautTemamChuntian}) that      there exists a unique solution to the parabolic problem which is sufficiently regular for all the subsequent calculations to be valid; in particular,
\begin{equation}\label{uep}
u^\epsilon \in   L^2 (0, T; H^4(\dom) ) \cap \mathcal {C}^1 ([0, T]; \Xi).
\end{equation}

\subsection{Global bounds independent of $\epsilon$  for $u^\epsilon$}
\label{secglobal}
Firstly, we prove the following global bounds:
\begin{lem}\label{main0}
We assume that
\begin{equation}\label{initial2}
u_{0} \in L^2 (\dom) ,
\end{equation}
\begin{equation}\label{f}
f \in L^2(0, T; L^2 (\dom)),
\end{equation}
then,  for every $T>0$ the following estimates independent of $\epsilon$ hold:
\begin{equation}\label{set1}
\begin{cases}
\quad u^\epsilon \,\,\mbox{is\,\, bounded\,\,in}\,\,L^\infty(0, T; \,L^2(\dom)),\\
\quad  u^\epsilon\mbox{ is\,\, bounded\,\,    in   }\,\,L^2(0, T; \,H^1 _0(\dom)),\\
\quad u^\epsilon _x (0, \cdot, \cdot ) \,\,\mbox{is\,\, bounded\,\,   in}\,\,L^2(0, T; \,L^2(I_{x^\perp})).
\end{cases}
\end{equation}
\end{lem}
\noindent \textbf{Proof.} As in \cite{SautTemamChuntian}, we multiply (\ref{10-11}) with $u^\epsilon$ and $xu^\epsilon$, integrate over $\dom$ and integrate by parts.
\qed\\\\

\subsection{Local bounds independent of $\epsilon$ for $u^\epsilon$}
\label{seclocal}

We first introduce a useful result:
\begin{lem}\label{localbounds1}
Under the same assumptions as in Lemma \ref{main0}, if we further suppose that
\begin{equation}\label{f1}
f \in L^\infty(0, T; L^2 (\dom)),
\end{equation}
 then we have
\begin{equation}\label{main8}
 | u^\epsilon_x(t) |^{2} \leq  |u^\epsilon_t(t)|^2 + \kappa ,\,\,\,\,\,\,\,\,\,0\leq t\leq T ,
\end{equation}
where
 $\kappa$ is a constant depending only on $|u_0|$, $|f|_{L^\infty (0, T; L^2 (\dom))}$ and $T$.
\end{lem}

\noindent\textbf{Proof.}
We rewrite (\ref{10-11})$_1$ as
\begin{eqnarray}\label{10-113}\displaystyle
\Delta \dfrac{\partial u^\epsilon }{\partial x}+   c\dfrac{\partial u^\epsilon }{\partial x}  +   u^\epsilon \dfrac{\partial u^\epsilon }{\partial x}  + \epsilon \, Lu^\epsilon      = - \dfrac{\partial u^\epsilon}{\partial t}+ f.
\end{eqnarray}
We multiply (\ref{10-113}) by $(1+x)u^\epsilon$, integrate over $\dom$, integrate by parts, and follow the same calculations as  in \cite{SautTemamChuntian}; we  find  when  $\epsilon \leq \frac{1}{4}$,
\begin{equation}\label{eq10-30}
\begin{split}
 |u_x^\epsilon |^{2}&
\leq - \int_{\dom} u^\epsilon_t \, (1+x)u^\epsilon\, d\dom  + |f|^2 + | (1+x)u^\epsilon|^2 + c |u^\epsilon|^2+     c^\prime |u^\epsilon |^6\\
 &\leq  |u^\epsilon_t|^2 + |f|^2  +    c^\prime |u^\epsilon|^2   +   c^\prime |u^\epsilon |^6,
 \end{split}
\end{equation}
Here and below $c^\prime$ indicates an absolute constant which may be different at each occurrence.
Hence if we call $\nu$ a bound of $|u^\epsilon|_{L^\infty(0, T; L^2 (\dom))} $ as in (\ref{set1})$_1$, we can set
\begin{equation}\label{kappa1}
\kappa =   |f|^2_{L^\infty(0, T; L^2 (\dom))}   +   c^\prime \nu ^2  +   c^\prime  \nu^6 + | u_{0x} |^{2};
\end{equation}
and  by (\ref{f1})  we obtain (\ref{main8}).
Thus we have completed the proof of Lemma \ref{localbounds1}.
\qed\\

Now we are ready to prove  the following result giving the local bounds on $u^\epsilon$ independent of $\epsilon$:
\begin{prop}\label{localbounds}
Under the same assumptions as in Lemma \ref{localbounds1},  if we further  suppose that
\begin{equation}\label{f2}
f_t \in L^\infty(0, T; L^2 (\dom)),
\end{equation}
\begin{equation}\label{f3}
f \in L^2(0, T; L^2 (I_x; H^2 (I_{x^\perp}))),
\end{equation}
\begin{equation}\label{initial1}
Lu_0 \in L^2 (\dom),
\end{equation}
\begin{equation}\label{initial3}
\Delta u_{0x} + u_0 u_{0x} + cu_{0x} - f(0)\in L^2 (\dom),
\end{equation}
\begin{equation}\label{uy0}
\nabla^\perp u_{0}, u_{0yy}, u_{0zz}\in L^2 (\dom),
\end{equation}
then there exists $T_{\ast}= \min (T, T_1) $,
\begin{equation}\label{tt}
T_1 = \dfrac{c_3}{\mu^4},
\end{equation}
\begin{equation*}
\mu= \mu (\kappa,    |f_t|_{L^\infty (0, T; L^2 (\dom))}, |Lu_0|,|\Delta u_{0x} + u_0 u_{0x} + cu_{0x} - f(0)|),
\end{equation*}
   such that for every $t$,  $0\leq t \leq T_{\ast}$,
\begin{equation}\label{v01}
 |u^\epsilon_t (t)|  \lesssim \mu,
\end{equation}
\begin{equation}\label{v02}
\int^{T_{\ast}}_0 |\nabla u^\epsilon_t (s)| ^2 \, ds \lesssim \mu , 
\end{equation}
\begin{equation}\label{mu1}
|\nabla u^\epsilon (t) | \leq  C(\mu),   
\end{equation}
\begin{equation}\label{y2}
\int^{T_{\ast}}_0 |\nabla u^\epsilon_y (s) |^2\,ds \leq   C(\mu),\, \,\,\, \,\, \int^{T_{\ast}}_0 |\nabla u^\epsilon_z (s) |^2\,ds \leq   C(\mu), 
\end{equation}
\begin{equation}\label{y3}
| u^\epsilon _{yy}(t) |\leq    C(\mu) , \,\, \,\,\,\,\,\,\,\,\,\,\,\,\,\,\,| u^\epsilon _{zz}(t) | \leq   C(\mu)  ,
\end{equation}
\begin{equation}\label{y4}
\int^{T_{\ast}}_0 |\nabla  u_{yy}^\epsilon (s) |^2\,ds \leq    C(\mu)  , \,\,\,\,\,\,\, \int^{T_{\ast}}_0 |\nabla  u_{zz}^\epsilon (s) |^2\,ds \leq  C(\mu), 
\end{equation}
\begin{equation}\label{yz1}
 \epsilon  \int^{T_{\ast}}_0[u^\epsilon_{yy}]_2^2\,\,ds \leq   C(\mu)  ,\,\,\,\,\,\,\,\,\,\epsilon \int^{T_{\ast}}_0 [u^\epsilon_{zz}]_2^2\,\,ds \leq  C(\mu) ,
\end{equation}
where $\lesssim$ means $\leq$ up to a multiplicative constant independent of $\epsilon$,  the constant $c_3$ depends only on the data, and the constant  $C(\mu)$ depends only on $\mu$ and the data and may be different at each occurrence.
\end{prop}

\vskip 1 mm

\noindent\textbf{Proof.}
We differentiate  (\ref{10-11}) in $t$,   write $u^\epsilon_t=v^\epsilon$ and we find:
\begin{eqnarray}\label{10-112}\displaystyle
\begin{cases}
&\,\,\, \dfrac{\partial v^\epsilon}{\partial t}+ \Delta \dfrac{\partial v^\epsilon }{\partial x} +  c\dfrac{\partial v^\epsilon }{\partial x}+ u^\epsilon \dfrac{\partial v^\epsilon }{\partial x}  +  v^\epsilon \dfrac{\partial u^\epsilon }{\partial x}  + \epsilon \, L v^\epsilon      =f_t,\\
&\,\,\,v^\epsilon (0)=u^\epsilon_{t0} = - \epsilon Lu_0 - \Delta u_{0x} - u_0 u_{0x}-cu_{0x} +f(0).
\end{cases}
\end{eqnarray}
Thus when $\epsilon \leq 1$,
\begin{equation}\label{v0}
|u^\epsilon_{t0} | \leq  |Lu_0| + | \Delta u_{0x} + u_0 u_{0x} + cu_{0x} - f(0)|.
\end{equation}
From (\ref{initial1}) and  (\ref{initial3}), we obtain
\begin{equation}\label{v011}
u^\epsilon_{t0} \mbox{ is bounded independently of $\epsilon$ in } L^2 (\dom) .
\end{equation}
Multiplying  (\ref{10-112}) by  $  (1+x)v^\epsilon$, integrating over $\dom$ and integrating by parts, dropping $\epsilon$ for the moment we find
\begin{equation*}
\begin{split}
  & \bullet\,\,  \int_{\dom} \frac{\partial v}{\partial t} (1+x)v \,d\dom = \dfrac{1}{2} \dfrac {d}{dt} |\sqrt{1+x}\, v|^2,\\
   & \bullet\,\,  \int_{\dom}  \Delta \dfrac{\partial v }{\partial x}   \, v \, d \dom =    \dfrac{1}{2} \int_{I_{x^\perp}} (v_{x}\big|_{x=0})^2\,dx^\perp ,\\
 & \bullet\,\,  \int_{\dom}  \Delta \dfrac{\partial v }{\partial x}   \, xv \, d\, \dom =  \frac{3}{2}|v _x|^2 +\frac{1}{2}| \nabla ^\perp v |^2, \\
  & \bullet \,\,\int_{\dom}  c v_x  (1+x)v       \,d\dom = - \frac{c}{2} \int_{\dom} v^2     \,d\dom,\\
  & \bullet\,\, \int_{\dom}    u v_x\,   (1+x)v       \,d\dom =  - \frac{1}{2} \int_{\dom} (1+x) u_x v^2      \,d\dom - \frac{1}{2}\int_{\dom} u v^2      \,d\dom,  \\
  & \bullet\,\, \int_{\dom}    v u_x\,   (1+x)v       \,d\dom =  \int_{\dom} (1+x) u_x v^2      \,d\dom,\\
  & \bullet\,\,  \int_{\dom} f_t (1+x)v \, d \dom \leq \dfrac{1}{2}  |f_t|^2 + \dfrac{1}{2}  |(1+x)v|^2 \leq \dfrac{1}{2}  |f_t|^2 +  |v|^2  ,\\
& \bullet\,\, \epsilon \int_{\dom} \dfrac{\partial ^4 v}{\partial x^4} (1+x) v \, dx \,dx^\perp=  \epsilon | \sqrt {1+x}\, v_{xx}  |^2 -
 \epsilon \int _{ I_{x^\perp}   }(v_x \big|_{x=0})^2 \,dx^\perp,\\
         & \bullet\,\,    \epsilon \int _{\dom} ( \dfrac{\partial^4 v}{\partial y^4} +
      \dfrac{\partial^4 v}{\partial z^4} )\, (1+x) v \,d \dom  =\epsilon \left( \left|\sqrt {1+x} \, v_{yy}\right|^2 +
\left|\sqrt {1+x}\, v_{zz}\right|^2 \right) .
\end{split}
\end{equation*}
Hence we arrive, when $\epsilon \leq \frac{1}{4}$, at
\begin{equation}\label{main2}
\begin{split}
\dfrac{d}{dt} | \sqrt {1+x} \,u^\epsilon_t | ^2 &+ |\nabla u^\epsilon_t |^{2}  + \dfrac{1}{4}|u^\epsilon_{t x}|_{x=0}|^2 _{L^2 (I_{x^\perp})} \\
 &+  2 \epsilon \left(  |\sqrt {1+x}\, u^\epsilon_{txx}   | ^2  +  |\sqrt {1+x}\, u^\epsilon_{tyy}   | ^2  +  |\sqrt {1+x}\, u^\epsilon_{tzz}   | ^2    \right)\\
&\,\,\,\,\,\,\,\,\,\,\,\,\leq \left | \int_{\dom} \left((1+x) u^\epsilon_x -  u^\epsilon \right)(u^\epsilon_t)^2   \,d\dom \right| + (c +2) |u^\epsilon_t|^2 + |f_t|^2 .
\end{split}\end{equation}
For the first term on the right-hand-side of (\ref{main2}), we have
\begin{equation}\label{nonlinear1}
\begin{split}
 \quad  \left | \int_{\dom} \left((1+x) u^\epsilon_x  - u^\epsilon \right)(u^\epsilon_t)^2   \,d\dom \right|&\lesssim  ({\mbox{with }\sigma ^\epsilon (t):= |u^\epsilon_x| + |u^\epsilon| })\\
 & \lesssim \sigma^\epsilon(t) |u^\epsilon_t|^2 _{L^4(\dom)}\\
 & \lesssim \mbox{(by $H^{3/4} \subset L^4$ in $3D$)}\\
 & \lesssim \sigma^\epsilon(t) |u^\epsilon_t|^{1/2} |\nabla u^\epsilon_t|^{3/2}\\
 &   \leq c^\prime (\sigma ^\epsilon (t) )^4 |u^\epsilon_t|^2 + \frac{1}{8} |\nabla u^\epsilon_t|^2\\ 
 &  \leq (\mbox{by } (\ref{main8}) )  \\
& \leq c^\prime \left (  |u^\epsilon_t| ^6+   \kappa^2    |u^\epsilon_t|^2 + |u^\epsilon|^4   |u^\epsilon_t|^2 \right)  + \frac{1}{8} |\nabla u^\epsilon_t|^2  .
\end{split} \end{equation}
Applying (\ref{nonlinear1}) to (\ref{main2}),   we obtain
\begin{equation}\label{main10}
\begin{split}
\dfrac{d}{dt} | \sqrt {1+x} \,u^\epsilon_t |^2 &+ \dfrac{7}{8}  |\nabla u^\epsilon_t |^{2}   + \dfrac{1}{4}|u^\epsilon_{t x}|_{x=0}|^2 _{L^2 (I_{x^\perp})}\\
&+   2 \epsilon  \left(  |\sqrt {1+x}\,u^\epsilon_{yxx}   | ^2  +  |\sqrt {1+x}\,u^\epsilon_{yyy}   | ^2  +  |\sqrt {1+x} \, u^\epsilon_{yzz}   | ^2    \right)\\
& \,\,\, \,\,\,\,\,\,\,\,\,\,\, \,\,\,\,\leq c_1 ( |u^\epsilon_t|^2 +1)^3 + |f_t|^2   \\
 & \,\,\, \,\,\,\,\,\,\,\,\,\,\, \,\,\,\,\leq  c_1 ( | \sqrt {1+x}\, u^\epsilon_t|^2 +1)^3 + |f_t|^2  .
\end{split}\end{equation}
{where} $c_1$ depends only on $\kappa$.
Setting $ | \sqrt {1+x}\, u^\epsilon_t|^2 +1 :=\mathcal {Y}^\epsilon$,  (\ref{main10}) implies that
\begin{equation}\label{key}
\dfrac{d}{dt} \mathcal {Y}^\epsilon  \leq c_2 (\mathcal {Y}^\epsilon)^3,
\end{equation}
with  $c_2 := c_1 + |f_t|^2_{L^\infty (0, T; L^2 (\dom))} $.
Thus
\begin{equation}\label{v01y}
 \mathcal {Y}^\epsilon(t)  \leq 2\mu_0^2,  \,\,\,\,\,\,\,0\leq t \leq \dfrac{3 }{8c_2 \mu_0^4 },
\end{equation}
where $\mu_0$ is a bound of  $\sqrt {\mathcal {Y}^\epsilon(0)}$  independent of $\epsilon$ as provided by (\ref{v011}).
Now (\ref{v01y}) implies that
\begin{equation}\label{v0111}
 |u^\epsilon_t (t)|  \lesssim  {\mu_0},\,\,\,\, \,\,\,0\leq t \leq  \dfrac{3 }{8c_2 \mu_0^4}.
\end{equation}
Then by 
  (\ref{v0111}) and (\ref{main8}) we deduce that
\begin{equation}\label{mu5}
| u^\epsilon_x (t) | \lesssim \mu  ,\,\,\,\,\,\,\,0\leq t \leq  T_{\ast}.
\end{equation}
with $\mu:=  {\mu_0} + \sqrt{\kappa}$, and $T_{\ast} = \min (T, T_1)$,
\begin{equation}\label{cp}
T_1= \frac{c_3}{ \mu^4} \leq \frac{3 }{8c_2 \mu_0^4 } .
\end{equation}
By (\ref{v0111}), (\ref{cp}) and (\ref{main10}) we obtain (\ref{v02}).

We  multiply  (\ref{10-11})  by  $(1+x)u^\epsilon_{yy}$, integrate over $\dom$ and integrate by parts, dropping $\epsilon$ for the moment, we find
\begin{align*}\displaystyle
& \bullet \int_{\dom} u_{t} \, (1+x) u_{yy} \, d\dom = -\dfrac{1}{2} \dfrac{d}{dt} | \sqrt {1+ x}\, u_y   |^2,\\
&  \bullet \int_{\dom} u_{xxx} \,(1+x) u_{yy}\, d\dom =  -\frac{3}{2} |u^2_{xy} | - \frac{1}{2} \left|u^2 _{xy} \big|_{x=0}\right|_{L^2 (I_y)},\\
&  \bullet \int_{\dom} u_{xyy} \,(1+x) u_{yy}\, d\dom =  -\frac{1}{2} |u^2_{yy} |, \\
& \bullet \int_{\dom} u_{xzz} \,(1+x) u_{yy}\, d\dom =  -\frac{1}{2} |u^2_{zy} | ,
\end{align*}
\begin{align*}
& \bullet \int_{\dom} c u_x (1+x) u_{yy} \, d\dom = \frac{c}{2} \int_{\dom} u^2_{y} \, d\dom,\\
&  \bullet \int_{\dom} u u_{x} \,(1+x) u_{yy}\, d\dom =  \dfrac{1}{2} \int_{\dom} u^2_y \left (  u- (1+x) u_x   \right) \, d\dom,\\
&\bullet    \epsilon \int _{\dom } u_{xxxx}\, (1+x) u_{yy}\, d\dom = -  \epsilon \int _{\dom } u_{xxxxy} (1+x) u_{y}\, d\dom   \\
 &\,\,\,\quad\quad\quad\quad\quad\quad\quad\quad\quad\quad\quad\,\,\,\,=  \epsilon  \int _{\dom } u_{xxxy}  u_{y}\, d \dom   + \epsilon \int _{d\dom } u_{xxxy} (1+x) u_{xy} \,d\dom    \\
&\,\quad\quad\quad\quad\quad\quad\quad\quad\quad\quad\quad\quad   =  - 2\epsilon \int _{\dom } u_{xxy}  u_{xy} d\dom  -\epsilon  \int _{\dom }(1+x) u^2_{xxy}  d\dom \\
& \,\quad\quad\quad\quad\quad\quad\quad\quad\quad\quad\quad\quad= \epsilon \int _{\dom } u^2_{xy}\big|_{x=0} \,d\dom    -  \epsilon \int _{\dom }(1+x) u^2_{xxy}\, d\dom , \\
&\bullet    \epsilon \int _{\dom } u_{zzzz}\, (1+x) u_{yy}\, d\dom = -\epsilon \int _{\dom }(1+x) u^2_{zzy}\, d\dom , \\
& \bullet   \epsilon \int _{\dom } u_{yyyy} \,(1+x) u_{yy}\, d\dom = -  \epsilon  \int _{\dom }  (1+x) u^\epsilon_{yyy}\, d\dom , \\
& \bullet \int _{\dom } f  (1+x)\, u_{yy} \, d\dom = - \int _{\dom } f_y \, (1+x) u_{y} \, d\dom \leq  \frac{1}{2} |(1+x) u_y| ^2 + \frac{1}{2} |f_y|^2.
\end{align*}
Hence when $\epsilon \leq \frac{1}{4}$, we have
\begin{equation}\label{main5}
\begin{split}
\dfrac{d}{dt} | \sqrt {1+ x} \, u^\epsilon_y   |^2 &+  | \nabla u^\epsilon_y|^2 + \dfrac{1}{4} \left|u^\epsilon _{xy} \big|_{x=0}\right|^2_{L^2 (I_y)} \\
&+2 \epsilon \left(  |\sqrt {1+x}  \, u^\epsilon_{yxx}   | ^2  +  |\sqrt {1+x} \, u^\epsilon_{yyy}   | ^2  +  |\sqrt {1+x}\,  u^\epsilon_{yzz}   | ^2    \right)\\
& \,\,\, \,\,\,\,\leq   \left  |\int_{\dom} (u^\epsilon_y)^2 \left ( (1+x) u^\epsilon_x-  u^\epsilon \right) \, d\dom \right| + (c + 2)|u^\epsilon_y|^2 + |f_y|^2.
\end{split} \end{equation}
For  the   first term on the right-hand-side of (\ref{main5}), we find
\begin{equation}\label{nonlinear12}
\begin{split}
  \left |\int_{\dom} (u^\epsilon_y)^2 \left ( (1+x) u^\epsilon_x - u^\epsilon \right) \, d \dom \right|  & \lesssim  ( \mbox{with }{\sigma ^\epsilon (t):= |u^\epsilon_x| + |u^\epsilon| })\\
    & \lesssim \sigma^\epsilon(t) |u_y^\epsilon|^2 _{L^4(\dom)}\\
    &  \leq c^\prime (\sigma ^\epsilon (t))^4 | u_{y}^\epsilon|^2 + \frac{1}{8} |\nabla u_{y}^\epsilon|^2 \\
 &\leq(\mbox {by (\ref{mu5})})\\
&\leq c^\prime  \left( \mu^4  +  |u^\epsilon|^4
\right) |u^\epsilon_{y}|^2 + \frac{1}{8} |\nabla u_{y}^\epsilon|^2   ,\,\,\,\,0\leq t \leq T_{\ast}.
\end{split}\end{equation}
Applying (\ref{nonlinear12}) to (\ref{main5}),   we find
\begin{equation}\label{main7}
\begin{split}
 \dfrac{d}{dt} | \sqrt {1+ x}\,u^\epsilon_y   |^2 & + \dfrac{7}{8} | \nabla u^\epsilon_y|^2 + \dfrac{1}{4} \left|u^\epsilon _{xy} \big|_{x=0}\right|^2_{L^2 (I_y)}\\
& + 2 \epsilon \left(  |\sqrt {1+x} \, u^\epsilon_{yxx}   | ^2  +  |\sqrt {1+x}\,  u^\epsilon_{yyy}   | ^2  +  |\sqrt {1+x}\,  u^\epsilon_{yzz}   | ^2    \right)\\
& \,\,\,\,\,\,\,\,\,\,\, \leq  {c^\prime}  \mu^4  |u^\epsilon_{y}|^2 + |f_y|^2  \\
& \,\,\,\,\,\,\,\,\,\,\,\leq c^\prime \mu^4 |\sqrt {1+ x}\,u^\epsilon_{y}|^2 + |f_y|^2   ,\,\,\,\,\,\,\,\,\,\,\,\,0\leq t \leq T_{\ast}.
\end{split}\end{equation}
We can then  close the  Gronwall inequality on the time interval  $(0, T_{\ast})$, and obtain
\begin{align*}
 |\sqrt {1+ x}\,u^\epsilon_y (t)  |^2 &\leq  {C}( \mu) \left( |\sqrt {1+ x}\,u_{0y}|^2 + \int ^{T_{\ast}}_0 |f_y(s)|^2\,ds        \right)\\
& \leq C(\mu, |u_{0y}|, |f_y|_{L^2 (0, T; L^2 (\dom))}) ,\,\,\,\,\,\,\,\,\,\,0\leq t \leq T_{\ast},
\end{align*}
which implies
\begin{equation}\label{uu2}
|u^\epsilon_y (t) | \leq  C({\mu} ) ,\,\,\,\,\,\,\,0\leq t \leq T_{\ast}.
\end{equation}
By (\ref{uu2}) and (\ref{main7}) we obtain
\begin{equation}\label{deltauy}
\int^{T_{\ast}}_0 |\nabla u^\epsilon_y (s) |^2\,ds \leq C(\mu).
\end{equation}

Similarly, we can obtain the same kind of estimates for  $u^\epsilon_z $, $\nabla u^\epsilon_{z}$, that is
\begin{equation}\label{uu22}
|u^\epsilon_z (t) | \leq  C({\mu} ) ,\,\,\,\,\,\,\,0\leq t \leq T_{\ast},
\end{equation}
\begin{equation}\label{deltauy2}
\int^{T_{\ast}}_0 |\nabla u^\epsilon_z (s) |^2\,ds \leq C(\mu).
\end{equation}
From (\ref{mu5}), (\ref{uu2}) and (\ref{uu22}) we obtain (\ref{mu1}).

We then  multiply   (\ref{10-11}) by $ (1+x) u^\epsilon_{yyyy} $,   integrate over $\dom$ and integrate by parts, to find
\begin{equation*}\displaystyle
\begin{split}
& \bullet \int_{\dom} u_{t} \, (1+x) u_{yyyy} \, d\dom = \dfrac{1}{2} \dfrac{d}{dt} |  \sqrt {1+ x} \, u_{yy}   |^2,\\
&  \bullet \int_{\dom} u_{xxx} \,(1+x) u_{yyyy}\, d\dom =  \frac{3}{2} |u_{xyy} |^2 + \frac{1}{2} \left|u _{xyy} \big|_{x=0}\right|^2_{L^2 (I_y)},\\
&  \bullet \int_{\dom} u_{xyy} \,(1+x) u_{yyyy}\, d\dom = \frac{1}{2}|u^2_{yyy} |^2, \\
&  \bullet \int_{\dom} u_{xzz} \,(1+x) u_{yyyy}\, d\dom = \frac{1}{2}|u^2_{zyy} |^2,\\
&  \bullet \int_{\dom} cu_x \,(1+x) u_{yyyy}\, d\dom = - \frac{c}{2}  |u_{yy}|^2,\\
&  \bullet \int_{\dom} u u_{x} \,(1+x) u_{yyyy}\, d\dom \\
 &\,\,\,\,\,\,\,\,= - \int_{\dom} u_{y} u_x(1+x) u_{yyy}   \,d \dom - \int_{\dom} u u_{xy}(1+x) u_{yyy}   \,d \dom\\
& \,\,\,\,\,\,\,\, = \int_{\dom} u^2_{yy} u_x(1+x)    \,d \dom + 2 \int_{\dom} u_{y} u_{xy}(1+x)u_{yy}  \,d \dom  +
 \int_{\dom} u u_{xyy}(1+x)u_{yy}  \,d \dom  \end{split}
 \end{equation*}
 \begin{equation*}
 \begin{split}
 & \,\,\,\,\,\,\,\,=  \int_{\dom} u^2_{yy} u_x(1+x)    \,d \dom  + \int_{\dom} \frac{\partial \left (u^2 _{y} \right)}{\partial x}(1+x)u_{yy}  \,d \dom + \frac{1}{2} \int_{\dom} u  \frac{\partial\left (u^2_{yy} \right )}{\partial x} (1+x) \,d \dom \\
 & \,\,\,\,\,\,\,\,= \frac{1}{2} \int_{\dom} u^2_{yy} u_x(1+x)    \,d \dom  -  \int_{\dom} u^2_y u_{yy}  \,d \dom -  \int_{\dom} u^2_{y} (1+x) u_{xyy}   \,d \dom - \frac{1}{2}  \int_{\dom} u u^2_{yy}  \,d \dom \\
 &\,\,\,\,\,\,\,\,  = \frac{1}{2} \int_{\dom}\left ( u_x(1+x)   -u \right )u^2_{yy}   \,d \dom  -  \int_{\dom} u^2_y u_{yy}  \,d \dom  -  \int_{\dom} u^2_{y} (1+x) u_{xyy}   \,d \dom \\
& \bullet   \epsilon \int _{\dom } u_{xxxx}\, (1+x) u_{yyyy}\, d\dom = -  \epsilon \int _{\dom } u_{xxxxy} (1+x) u_{yyy}\, d\dom    \\
&\,\,\,\quad\quad\quad\quad\quad\quad\quad\quad\quad\quad\quad\quad\,\,=  - \epsilon \int _{\dom } u_{xxxy}  u_{yyy}\,d \dom   +  \epsilon\int _{\dom } u_{xxxy} (1+x) u_{xyyy} \, d\dom    \\
&\,\,\,\quad\quad\quad\quad\quad\quad\quad\quad\quad\quad\quad\quad\,\,= -  \epsilon \int _{I_{x^\perp} } u^2_{xyy}\big|_{x=0} \, d\, I_{x^\perp}    + \epsilon \int _{\dom }(1+x) u^2_{xxyy}\, d\dom , \\
& \bullet    \epsilon \int _{\dom } u_{yyyy}\, (1+x) u_{yyyy}\, d \dom = \epsilon |\sqrt{1+x}\, u_{yyyy}|^2 , \\
& \bullet    \epsilon \int _{\dom } u_{zzzz} \, (1+x) u_{yyyy}\, d \dom = \epsilon |\sqrt{1+x} \,u_{zzyy}|^2 , \\
& \bullet \int _{\dom } f  (1+x) u_{yyyy} \, d\dom = - \int _{\dom } f_{yy} \, (1+x) u_{yy} \, d\dom \leq  \frac{1}{2} |(1+x) u_{yy}| ^2 + \frac{1}{2} |f_{yy}|^2.
\end{split}
\end{equation*}
Hence  when $\epsilon \leq \frac{1}{4}$,
\begin{equation}\label{main1}
\begin{split}
\dfrac{d}{dt} |\sqrt {1+ x}&\,u^\epsilon_{yy}   |^2 +  | \nabla u^\epsilon_{yy}|^2 + \dfrac{1}{4} |u^\epsilon _{xyy} \big|^2_{x=0}|_{L^2 (I_y)} \\
&+ 2 \epsilon \left(  |\sqrt {1+x}  \,u^\epsilon_{yyxx}   | ^2  +  |\sqrt {1+x}\,  u^\epsilon_{yyyy}   | ^2  +  |\sqrt {1+x} \,u^\epsilon_{yyzz}   | ^2    \right)\\
&\,\,\,\,\,\,\,\,\,\,\,\,\,\leq   \left | \int_{\dom} \left ( (1+x) u^\epsilon_x - u^\epsilon\right) (u^\epsilon_{yy})^2\, d\dom \right|  +     2 \left |\int_{\dom} u^2_y u_{yy}  \,d \dom \right | \\
 & \,\,\,\,\,\,\,\,\,\,\,\,\,\,\,\,\,\,\,\,\,\,\,\,\,+ 2\left | \int_{\dom} u^2_{y} (1+x) u_{xyy}   \,d \dom \right | + (c +1) |u^\epsilon_{yy}|^2 + |f_{yy}|^2 \\
 & \,\,\,\,\,\,\,\,\,\,\,\,\,\,:= I^\epsilon_1 + I^\epsilon_2 + I^\epsilon_3 + (c +1) |u^\epsilon_{yy}|^2 + |f_{yy}|^2.
\end{split} \end{equation}
For $I^\epsilon_1$, by the similar calculations in (\ref{nonlinear12}) we deduce
\begin{equation}\label{I1}
 I^\epsilon_1 \leq c^\prime \left( \mu^4  +  |u^\epsilon|^4
\right) |u^\epsilon_{yy}|^2 + \frac{1}{8} |\nabla u_{yy}^\epsilon|^2   ,\,\,\,\,0\leq t \leq T_{\ast}.
\end{equation}
For $I^\epsilon_2$ we have
\begin{equation}\label{I2}
\begin{split}
  I^\epsilon_2 &\leq 2 | u^\epsilon_y|^2 _{L^4 (\dom)} |u^\epsilon _{yy}|\\
  & \lesssim | u^\epsilon_y|^{1/2}  | \nabla u^\epsilon_y|^{3/2}       |u^\epsilon _{yy}|\\
  & \leq ( \mbox{by  }(\ref{uu2}) )\\
  & \leq C(\mu)   | \nabla u^\epsilon_y|^{3/2}       |u^\epsilon _{yy}|\\
  & \leq C(\mu)   | \nabla u^\epsilon_y|^{3/2}   |u^\epsilon _{yy}| ^2 + C(\mu )   | \nabla u^\epsilon_y|^{3/2} , \,\,\,\,0\leq t \leq T_{\ast}.
  \end{split}\end{equation}
For $I^\epsilon_3$ we have
\begin{equation}\label{I3}
\begin{split}
  I^\epsilon_3 &\leq 2 | u^\epsilon_y|^2 _{L^4 (\dom)} |u^\epsilon _{xyy}|\\
  & \lesssim | u^\epsilon_y|^{1/2}  | \nabla u^\epsilon_y|^{3/2}       |u^\epsilon _{xyy}|\\
  & \leq \frac{1}{8}  |u^\epsilon _{xyy}|^2 + c ^\prime   | u^\epsilon_y|  | \nabla u^\epsilon_y|^{3}    \\
  & \leq (\mbox{by Jensen's inequaliy, }  | \nabla u^\epsilon_y|^{3} \lesssim | u^\epsilon_{xy}   | ^ 3 + | u^\epsilon_{yy} |^3 + | u^\epsilon_{yz}    |^3  )\\
  & \leq \frac{1}{8}  |u^\epsilon _{xyy}|^2 + c ^\prime   | u^\epsilon_y|  |u^\epsilon_{xy}|^{3} + c ^\prime   | u^\epsilon_y|  |u^\epsilon_{yy}|^{3} + c ^\prime   | u^\epsilon_y|  |u^\epsilon_{yz}|^{3}  \\
  &:= \frac{1}{8}  |u^\epsilon _{xyy}|^2 + J^\epsilon_4 + J^\epsilon_5 + J^\epsilon_6.
  \end{split}\end{equation}
We now estimate $J^\epsilon_4$.  {We observe that} since  $u^\epsilon _x =0$ at $y = \pm \frac{\pi}{2}$,
\begin{equation}\label{yel}
\begin{split}
|u^\epsilon_{xy}|^2  = \int _{\dom} \left (u^\epsilon_{xy}  \right )^2\, d \dom
 = - \int_{\dom }u^\epsilon_{x}u^\epsilon_{xyy}  \, d \dom
\leq c^\prime |u^\epsilon_{x}  | |u^\epsilon_{xyy}  | .
\end{split}
\end{equation}
Thus we have
\begin{equation}\label{j4}
\begin{split}
J^\epsilon_4 &\leq c^\prime   | u^\epsilon_y|    |u^\epsilon _{x}|^{3/2}  |u^\epsilon _{xyy}|^{3/2} \\
& \leq c^\prime |\nabla u^\epsilon |^{5/2} |u^\epsilon _{xyy}|^{3/2} \\
 &\leq (\mbox{by (\ref{mu1}) which is already proven)})\\
 & \leq C(\mu)^{5/2} |u^\epsilon _{xyy}|^{3/2}\\
 &\leq   C(\mu)^{10}   + \frac{1}{8} |u^\epsilon _{xyy}|^{2}, \,\,\,\,0\leq t \leq T_{\ast}.
  \end{split}\end{equation}

Similarly for $J^\epsilon_6$, since $u^\epsilon _z =0$ at $y = \pm \frac{\pi}{2} $, we can apply the intermediate derivative theorem to $u^\epsilon_z$,  and deduce that $|u^\epsilon_{zy}|^2 \leq c^\prime |u^\epsilon_{z}  | |u^\epsilon_{zyy}  |$. Hence by estimates similar as in (\ref{j4}) we have
\begin{equation}\label{j6}
J^\epsilon_6 \leq  C(\mu)^{10}   + \frac{1}{8} |u^\epsilon _{zyy}|^{2}, \,\,\,\,0\leq t \leq T_{\ast}.
\end{equation}

To estimate $J^\epsilon_5$, by (\ref{uu2}) we have
\begin{equation}\label{j5}
J^\epsilon_5 \leq   C(\mu) |u^\epsilon _{yy}|^{3}, \,\,\,\,\,\,\,\,\,0\leq t \leq T_{\ast}.
\end{equation}

Collecting the estimates in (\ref{j4}), (\ref{j5}) and (\ref{j6}), along with (\ref{I3}) we obtain
\begin{equation}\label{I33}
\begin{split}
  I^\epsilon_3 \leq  \frac{3}{8}  |\nabla u^\epsilon _{yy}|^2 + C(\mu)  +   C(\mu) |u^\epsilon _{yy}|^{3}
  , \,\,\,\,\,\,\,\,\,0\leq t \leq T_{\ast}.
  \end{split}\end{equation}

Collecting the estimates in (\ref{I1}), (\ref{I2}) and (\ref{I33}), along with (\ref{main1}) we obtain
\begin{equation}\label{main111}
\begin{split}
\dfrac{d}{dt} | \sqrt {1+ x}\, u^\epsilon_{yy}   |^2 &+ \frac{1}{2} | \nabla u^\epsilon_{yy}|^2 + 2 \epsilon \left(  |\sqrt {1+x}  \,u^\epsilon_{yyxx}   | ^2  +  |\sqrt {1+x}\,  u^\epsilon_{yyyy}   | ^2  +  |\sqrt {1+x} \,u^\epsilon_{yyzz}   | ^2    \right)\\
&\,\,\,\,\,\,\,\,\,\,\,\,\leq   c^\prime \left( \mu^4  +  |u^\epsilon|^4 + C(\mu)   | \nabla u^\epsilon_y|^{3/2} + C(\mu) |u^\epsilon _{yy}| + c +1 \right) |u^\epsilon_{yy}|^2\\
& \,\,\,\,\,\,\,\,\,\,\,\,\,\,\,\,\,\,\,\,\,\,\,\,+   C(\mu )   | \nabla u^\epsilon_y|^{3/2} +   C(\mu)  + |f_{yy}|^2,\,\,\,\,\,\,\,\,\,\,\,\,\,\,0\leq t \leq T_{\ast}.
\end{split} \end{equation}
In particular, setting $\eta ^\epsilon (t) =  c^\prime \left( \mu^4  +  |u^\epsilon|^4 + C(\mu)   | \nabla u^\epsilon_y|^{3/2} + C(\mu) |u^\epsilon _{yy}| + c +1 \right)$, from (\ref{main111}) we infer that
 \begin{equation}\label{main12}
\begin{split}
\dfrac{d}{dt} | \sqrt {1+ x}\, u^\epsilon_{yy}   |^2 & \leq   \eta^\epsilon (t) | \sqrt {1+ x}\,  u^\epsilon_{yy}|^2+   C(\mu )   | \nabla u^\epsilon_y|^{3/2} +   C(\mu)  + |f_{yy}|^2,\,\,\,\,0\leq t \leq T_{\ast}.
\end{split} \end{equation}
Since $| \nabla u^\epsilon_y|^{3/2} \leq | \nabla u^\epsilon_y|^{2} + c^\prime    $, along with (\ref{deltauy}) we deduce
\begin{equation*}
\int^{T_{\ast}}_0  \eta^\epsilon(s)  \,ds \leq C(\mu) .
\end{equation*}
We can then  close the  Gronwall inequality on the time interval  $(0, T_{\ast})$ in (\ref{main111}), and obtain
\begin{equation*}
|\sqrt {1+x}\,u^\epsilon_{yy} (t) | \leq   C(\mu, |u_{0yy}|, |f_{yy}|_{L^2 (0, T; L^2 (\dom))}) ,\,\,\,\,\,\,\,0\leq t \leq T_{\ast},
\end{equation*}
which implies
\begin{equation}\label{uu1}
|u^\epsilon_{yy} (t) | \leq   C(\mu) ,\,\,\,\,\,\,\,0\leq t \leq T_{\ast}.
\end{equation}
By (\ref{uu1}) and (\ref{main1}) we obtain
\begin{equation}\label{y8}
\int^{T_{\ast}}_0 |\nabla u^\epsilon_{yy} (s) |^2\,ds  \leq C(\mu) ,
\end{equation}
\begin{equation}\label{yz2}
 \epsilon  \int^{T_{\ast}}_0[u^\epsilon_{yy}]_2^2\,\,ds \leq C(\mu) ,\,\,\,\,\,\,\,0\leq t \leq T_{\ast}.
\end{equation}

Similarly we can obtain the same kind of estimates for  $u^\epsilon_{zz}$, $\nabla u^\epsilon_{zz} $ and $\epsilon  [u^\epsilon_{zz}]_2^2$.

Combining all the previous local bounds,  we obtain (\ref{v01})-(\ref{yz1}).
Hence we have completed the proof of Proposition \ref{localbounds}. \qed

%
%
%
%

\vskip 5 mm

\subsection{A singular perturbation argument}
\label{sing}

We are now ready to show the local estimates for $u^\epsilon_{xx}$ and $u^\epsilon\,u^\epsilon_x$ by singular perturbation.
\begin{prop}\label{localbound2}
Under the same assumptions as in Proposition \ref{localbounds}, we have
\begin{equation}\label{h2}
u_{xx}^\epsilon \mbox{ is bounded independently of $\epsilon$ in } L^2 (0, T_{\ast},\,L^2 (\dom)),
\end{equation}
\begin{equation}\label{uux}
u^\epsilon u^\epsilon_x \mbox{ is bounded independently of $\epsilon$    in } L^2 (0, T_{\ast}; L^2 (\dom)) .
\end{equation}


\end{prop}

\vskip 3 mm

\begin{rem}\label{rem1}
Note that  by (\ref{h2}) and (\ref{y3}) we deduce that
\begin{equation}\label{h22}
u^\epsilon \mbox{ is bounded independently of $\epsilon$ in } L^2 (0, T_{\ast},\,\Xi).
\end{equation}
\end{rem}

\begin{rem}
We know that
\begin{equation*}
\begin{split}
\int_{\dom}( u^\epsilon u^\epsilon_x)^{3/2} \, d\dom &\leq  \left (\int_{\dom} (u^\epsilon)^6 \, d \dom \right)^{1/4}  \left(\int_{\dom} (u^\epsilon_x)^2\,d\dom \right)^{3/4}\\
& = |u^\epsilon|^{3/2}_{L^6(\dom)} |u^\epsilon_x|^{3/2}\\
&\leq (\mbox{by $H^1(\dom) \subset L^6 (\dom)$ in $3D$})\\
 & \lesssim |\nabla u^\epsilon|^3 .
\end{split}
\end{equation*}
Hence
\begin{equation*}
 \esssup_{t \in (0, T_{\ast})} |u^\epsilon u^\epsilon_x(t)|_{L^{3/2}(\dom)} \lesssim \esssup _{t \in (0, T_{\ast})} |\nabla u^\epsilon (t)|^2 \lesssim (\mbox{by (\ref{mu1})}) \lesssim C(\mu)^{2},
\end{equation*}
which implies that
\begin{equation*}
\begin{split}
u^\epsilon u^\epsilon_x &\mbox{ is bounded independently of $\epsilon$    in } L^{\infty} (0, T_{\ast}; L^{3/2} (\dom)) ,\\
&\mbox{\,\,\,\, and hence    in } L^{3/2} (I_x; L^{3/2} ((0, T_{\ast}) \times I_{x^\perp})) .
\end{split}
\end{equation*}
Thus we can apply  Lemma \ref{A.2} in the Appendix with  $p=3/2$ and $Y=  L^{3/2} ((0, T_{\ast}) \times I_{x^\perp})$, and obtain
\begin{equation}\label{h21}
u_{xx}^\epsilon \mbox{ is bounded independently of $\epsilon$ in } L^{\infty}  (I_x; L^{3/2} ((0, T_{\ast}) \times I_{x^\perp})).
\end{equation}

However, to obtain more useful estimates as in (\ref{h2}) and (\ref{uux}), we need to use the following proof which provides a stronger result.
\end{rem}

\noindent\textbf{Proof of Proposition \ref{localbound2}.}
 We rewrite  the regularized equation  (\ref{10-11}) as  follows:
\begin{eqnarray}\label{x}\displaystyle
\begin{cases}
& u^\epsilon_{xxx} +  u^\epsilon u^\epsilon_x+ \epsilon  u^\epsilon_{xxxx} =  g^\epsilon , \\
& u^\epsilon(0)= u^\epsilon(1) = u^\epsilon_{x}(1) = u^\epsilon_{xx}(0)=0,
\end{cases}
\end{eqnarray}
where $g^\epsilon:= -u^\epsilon_t -  \Delta^\perp u^\epsilon_{x} - cu^\epsilon _x - \epsilon u^\epsilon_{yyyy} -  \epsilon u^\epsilon_{zzzz} +f $.
Hence by  (\ref{v01}), (\ref{y4}) and    (\ref{yz1}),    we know that each term in $g^\epsilon$ is bounded independently of $\epsilon$ in $ L^2 (0, T_{\ast},\,L^2 (\dom))$, and thus
\begin{equation}\label{g}
g^\epsilon \mbox{ is bounded independently of $\epsilon$  in } L^2 (0, T_{\ast},\,L^2 (\dom)) .
\end{equation}

Multiplying (\ref{x}) by $x$ and integrating in $x$ from $0$ to $1$, we find
\begin{align*}
& \bullet \int ^1 _ {0} x u_{xxx}\, dx = - \int ^1 _ {0}     u_{xx} \, dx +  u_{xx} x \big|_{x=0} ^ {x=1}     = u_x\big|_{x=0}  +  u_{xx}\big|_{x=1}  ,    \\
&  \bullet  \int^1_ {0  } x  u u_x  \, dx =\int^1_ {0  } \frac{\partial }{\partial x}\left( \frac {u^2}{2} \right) dx   =     -\dfrac {1}{2} \int^1_ {0  }  u^2\, dx  ,\\
&  \bullet   \epsilon \int ^1 _  {0   } x u_{xxxx}\, dx  =    - \epsilon \int^1_ {0  }u_{xxx}\,dx +   \epsilon  u_{xxx} x  \big|^{x=1}_{x=0} =  -  \epsilon u_{xx}\big|_{x=1}  +  \epsilon u_{xxx} \big|_{x=1}.
\end{align*}
Hence
\begin{equation}\label{eq2}
u^\epsilon_x\big|_{x=0}  +  u^\epsilon_{xx}\big|_{x=1}     -\dfrac {1}{2} \int^1_ {0  }  (u^\epsilon)^2\, dx   - \epsilon u^\epsilon_{xx}\big|_{x=1}  + \epsilon u^\epsilon_{xxx} \big|_{x=1} =  \int ^1 _ {0}g^\epsilon x \, dx.
\end{equation}

Integrating (\ref{x}) in $x$ from $\tilde x $ to $1$,   we obtain
\begin{align*}
& \bullet \int ^1 _ {\tilde x  } u_{xxx}\, dx =u_{xx} \big|_{x = 1} - u_{xx}    ,    \\
&  \bullet  \int^1_ {\tilde x  }  u u_x  \, dx =    -\dfrac {1}{2}u^2 ,\\
&  \bullet   \epsilon  \int ^1 _  {\tilde x  }  u_{xxxx}\, dx =  \epsilon u_{xxx}  \big|_{x=1} -  \epsilon u_{xxx}.
\end{align*}
Hence
\begin{equation}\label{a1}
u^\epsilon_{xx} \big|_{x = 1} - u^\epsilon_{xx}   -\dfrac {1}{2}(u^\epsilon)^2  +  \epsilon u^\epsilon_{xxx}  \big|_{x=1}   -  \epsilon u^\epsilon_{xxx} = \int^1 _{\tilde x} g^\epsilon \,dx.
\end{equation}
Then  (\ref{eq2})  and (\ref{a1}) imply
\begin{equation*}
u^\epsilon_x\big|_{x=0} -\dfrac {1}{2} \int^1_ {0  }  (u^\epsilon)^2\, dx - \epsilon u^\epsilon_{xx}\big|_{x=1} + u^\epsilon_{xx} + \dfrac {1}{2}(u^\epsilon)^2  + \epsilon u^\epsilon_{xxx}  =     \int ^1 _ {0} g^\epsilon x \, dx -   \int^1 _{\tilde x} g^\epsilon \,dx,
\end{equation*}
which we rewrite as
\begin{equation}\label{a2}
u^\epsilon_{xx} + \epsilon u^\epsilon_{xxx} = \epsilon u^\epsilon_{xx}\big|_{x=1} + h^\epsilon,
\end{equation}
where
\begin{equation}\label{hh}
h^\epsilon =   - u^\epsilon_x\big|_{x=0} +\dfrac {1}{2} \int^1_ {0  }  (u^\epsilon)^2\, dx   -\dfrac {1}{2}(u^\epsilon)^2       +     \int ^1 _ {0} g^\epsilon x \, dx -   \int^1 _{\tilde x} g^\epsilon \,dx.
\end{equation}
 Now we estimate the term $(u^\epsilon)^2  $ in (\ref{hh}). Since
  \begin{align*}
  |(u^\epsilon)^2|^2 \leq |u^\epsilon|^4_{L^4 (\dom)} \lesssim (\mbox{by $H^{3/4}(\dom) \subset L^4 (\dom)$ in 3D}) \lesssim  |\nabla u^\epsilon|^3 |u^\epsilon| ,
   \end{align*}
   we have
\begin{equation}
\int^ {T_{\ast}}_0|(u^\epsilon)^2|^2 \, ds   \lesssim \int^ {T_{\ast}}_0 |\nabla u^\epsilon|^3 |u^\epsilon|    \, ds \leq (\mbox{by (\ref{mu1}) and the Poincar\'e inequality} )\lesssim C(\mu)^4 T_{\ast}.
\end{equation}
Thus
\begin{align}\label{u2}
(u^ \epsilon)^2 \mbox{ is bounded independently of $\epsilon$   in } L^2 (0, T_{\ast}; L^2 (\dom)).
\end{align}
Applying   (\ref{set1})$_4$, (\ref{g}) and (\ref{u2}) to (\ref{hh}) we find
\begin{equation}\label{h}
h^\epsilon \mbox{ is bounded independently of $\epsilon$    in } L^2 (0, T_{\ast},\,L^2 (\dom)) .
\end{equation}

Multiplying (\ref{a2}) by $u^\epsilon_{xx}$, integrating in $x$ from $0$ to $1$, we obtain
\begin{align*}
& \bullet \epsilon \int ^1 _ {0}  u_{xxx} u_{xx}  \, dx = \frac{\epsilon}{2}  u^2_{xx}\big|_{x=1}  ,    \\
&  \bullet \epsilon \int^1_ {0  }u_{xx}\big|_{x=1}\, u_{xx} \, dx =   \epsilon u_{xx}\big|_{x=1}\, \int^1_ {0  } u_{xx} \, dx     = - \epsilon u_{xx}\big|_{x=1} \, u_x \big|_{x=0};
\end{align*}
hence we arrive at
\begin{align*}
\int ^1 _0 (u^\epsilon_{xx})^2 \, dx + \frac{\epsilon}{2} ( u^\epsilon_{xx} )^2 \big|_{x=1}& = - \epsilon u^\epsilon_{xx} \big|_{x=1} u^\epsilon_x \big|_{x=0}+ \int^1 _0 u^\epsilon_{xx}\, h^\epsilon \, dx,\\
& \leq \frac{\epsilon}{4} ( u^\epsilon_{xx}\big|_{x=1} )^2 + c^{\prime } \epsilon (u^\epsilon_x \big|_{x=0})^2 + \frac{1}{2} |u^\epsilon_{xx}|^2_{L^2 (I_x)} + \frac{1}{2}  |h^\epsilon|^2_{L^2 (I_x)}.
\end{align*}
%
Thus
\begin{equation}\label{uxx2}
\frac{1}{2}\int ^1 _0 (u^\epsilon_{xx})^2 \, dx + \frac{\epsilon}{4} ( u^\epsilon_{xx} \big|_{x=1}  )^2 \leq  c^{\prime} \epsilon (u^\epsilon_x \big|_{x=0})^2  + \frac{1}{2}  |h^\epsilon|^2_{L^2 (I_x)}.
\end{equation}
We integrate both sides of (\ref{uxx2}) in  $I_{x^\perp}$ and then in time from $0$ to $T_{\ast}$; by   (\ref{h}) and (\ref{set1})$_4$ we obtain (\ref{h2}). As in Remark \ref{rem1}, we thus have (\ref{h22}).

Now since
\begin{align*}
|u^\epsilon u^\epsilon_x|^2 & \leq   |u ^\epsilon|^2_{L^4 (\dom)}|u^\epsilon_x|^2_{L^4 (\dom)}\\
&\lesssim (\mbox{by $H^{3/4}(\dom) \subset L^4 (\dom)$ in 3D}) \\
&\lesssim  |u^\epsilon|^{1/2} |\nabla u^\epsilon|^{3/2}  |u_x^\epsilon| ^{1/2} [u^\epsilon]_2^{3/2}\\
 &\lesssim |u^\epsilon|^2 |\nabla u^\epsilon| ^4   + [u^\epsilon]_2 ^2,
\end{align*}
hence we obtain
\begin{align*}
\int^{T_{\ast}} _0|u^\epsilon u^\epsilon_x|^2\, ds  &\lesssim \int^{T_{\ast}}_0 |u^\epsilon|^2 |\nabla u^\epsilon| ^4\, ds  +  \int^{T_{\ast}}_0 [u]_2 ^2\, ds.
\end{align*}
This together with  (\ref{mu1}) and (\ref{h22}) implies (\ref{uux}).
\qed

%
%
%
%

\section{Passage to the limit}
\label{secexistence}

%
%

Using a compactness argument, we can pass to the limit in  (\ref{10-11}). Hence we obtain (\ref{eq0}), with a function $u\in  {\mathcal C}^1 ([0, T_{\ast}]; L^2 (\dom) ) \cap  L^2 (0, T;  H^1_0(\dom))$.
Then we rewrite  (\ref{eq0})  as
\begin{equation}\label{eq00}
u_{xxx} = -   u_t - \Delta ^{\perp} u_x -  c u_x - u u_x  -f.
\end{equation}
From (\ref{v01}), (\ref{y4}) and (\ref{uux}), we infer that each term in the right-hand-side of (\ref{eq00}) belongs to  $ L^2 (0, T_{\ast}; L^2 (\dom))$, and hence
\begin{equation}\label{uxxx}
u_{xxx} \in L^2 (0, T_{\ast}; L^2 (\dom)).
\end{equation}

Now we are ready to state the main result: the local existence of  strong solutions.
\begin{thm}\label{localstrong}
 The assumptions are the same as in Proposition \ref{localbounds}, that is (\ref{initial2}), (\ref{f}), (\ref{f1}) and  (\ref{f2})-(\ref{uy0}). We suppose also that the following compatibility conditions hold:
\begin{equation}\label{compatibility}
u_0  = 0 \mbox{ on } \partial \dom,\,\,u_{0x} \big|_{x=1} =0,\,\,  u_{0yy} \big|_{y= \pm  \frac{\pi}{2}} =u_{0zz}  \big|_{z= \pm \frac{\pi}{2}} =0,
\end{equation}
\begin{equation}\label{compatibility1}
u_{t0}= 0 \mbox{ on } \partial \dom,\,\,\frac{\partial u_{t0} }{\partial x} \Big|_{x=1} =0,\,\, \frac{\partial^2 u_{t0} }{\partial y^2}  \Big|_{y= \pm  \frac{\pi}{2}} =\frac{\partial^2 u_{t0} }{\partial z^2}   \Big|_{z= \pm \frac{\pi}{2}} =0,
\end{equation}
where $u_{t0}=- \Delta u_{0x} - u_0 u_{0x}-cu_{0x} +f(0) $.
Then there exists a local strong solution  to  (\ref{eq0})-(\ref{eq107})  on some time interval $[0, T_{\ast})$,  $T_{\ast} >0$ depending only on the data as in Proposition \ref{localbounds},  such that
\begin{equation}\label{defs}
 \nabla u, u_{yy}, u_{zz}, u_t  \in  L^\infty(0, T_{\ast}; L^2 (\dom)) ,
\end{equation}
\begin{equation}\label{da}
u \in L^2 (0, T_{\ast}; D(A)),
\end{equation}
\begin{equation}\label{da11}
u \in L^2 (0, T_{\ast}; \Xi),
\end{equation}
\begin{equation}\label{defs2}
u \in L^2 (0, T_{\ast};   H^3 (I_x; L^2 (I_{x^\perp})) \cap H^3 (I_{x^\perp}; L^2 (I_{x}))),
\end{equation}
\begin{equation}\label{def3}
u_t \in   L^2(0, T_{\ast}; H^1(\dom)) .
\end{equation}
Moreover, we have for every $t \in (0, T_{\ast})$,
\begin{equation}\label{uyy=0}
u_{yy}(t)  \big|_{y= \pm  \frac{\pi}{2}} =u_{zz}(t)  \big|_{z= \pm \frac{\pi}{2}} =0. 
\end{equation}

\end{thm}

\begin{rem}
We have proven that all the spatial derivatives of the third order of $u$ are in  $L^2(0, T_{\ast}; L^2 (\dom)) $, except for $u_{xxy}$ and $u_{xxz}$.
\end{rem}

\noindent\textbf{Proof.} We rewrite (\ref{eq0}) as
\begin{equation}\label{eq000}
Au = -   u_t -  u u_x  -f;
\end{equation}
 from (\ref{v01}) and  (\ref{uux}) we know that each term on the right-hand side of (\ref{eq000}) belongs to $L^2 (0,  T_{\ast}; L^2 (\dom))$. Hence $Au$ belongs to the same space. We also know that $u_x( 1, \, x^\perp,\, t )=0$, $t\in [0, T]$ using the same argument as in  \cite{SautTemamChuntian}. Hence we obtain   (\ref{da}). Now by (\ref{da2}) and (\ref{da}), we deduce (\ref{da11}).

By (\ref{y2}), we know that $u_{yyy}$, $u_{zzz}$ both belong to $L^2 (0, T_{\ast}; L^2 (\dom))$. Hence  we can apply  the trace theorem and pass to the limit on the boundary conditions in (\ref{10-4}) to obtain (\ref{uyy=0}).

The other results can be deduced directly from (\ref{v01})-(\ref{y4}) and (\ref{uxxx}).
\qed\\

%

\begin{rem}
As for the periodic case, that is, (\ref{eq0}) and the boundary and initial conditions (\ref{eq3}), (\ref{eq33}) and  (\ref{eq1-26p}), the results are exactly the same as in the Dirichlet case  discussed above. The reasoning is totally the same and therefore we skip it.
\end{rem}

\section{Appendix: a trace result}
\label{seca}
We recall a trace result from \cite{SautTemamChuntian}, which is used in the article.
\begin{lem}\label{A.2}
Let $Y$ be a {(not necessarily reflexive) Banach}    space and let ${p\geq 1}$. Assume that two sequences of functions $u^\epsilon$, $g^\epsilon \in L^p_x(I_x; Y)  $ satisfy
\begin{equation}\label{traceconv}
\begin{cases}
\quad u^\epsilon _{xxx} +\epsilon u^\epsilon_{xxxx}=g^\epsilon,\\
\quad u^\epsilon (0)= u^\epsilon (1)= u_x^\epsilon (1)= u_{xx}^\epsilon (0)=0,
\end{cases}
\end{equation}
 with $g^\epsilon$ bounded in $L^p_x(I_x; Y)  $  as $\epsilon \rightarrow 0$. Then $u^\epsilon _{xx}$ (and hence $u^\epsilon_x$, and $u^\epsilon$) is bounded in $L^\infty _x(I_x; Y)  $ as $\epsilon \rightarrow 0$. {Furthermore if Y is reflexive}, then for any subsequence $u^ \epsilon \rightarrow u$ converging {strongly or weakly in $L^q_x(I_x; Y) $, $  1 < q < \infty$}, $u^\epsilon_x (1)$ converges to $u_x (1)$ in $Y$ (weakly at least), and hence $u_x (1)=0$.
\end{lem}

\section*{Acknowledgments}
This work was partially supported by the National Science Foundation under the grants, DMS-0906440  and DMS 1206438, and by the Research Fund of Indiana University.

The author would like to thank Professor Roger Temam, Professor Jean-Claude Saut, and
Professor Nathan Glatt-Holtz for their encouragements and suggestions.

   \newpage

\footnotesize
\bibliographystyle{amsalpha}

\bibliography{ref-3}

\normalsize

%

\end{document}